\documentclass[reqno]{amsart}

\usepackage{amsmath}
\usepackage{amssymb,latexsym,amstext,amsthm}
\usepackage{verbatim} 
\usepackage{eucal} 

\theoremstyle{plain}
\newtheorem{thm}[equation]{Theorem}
\newtheorem{pro}[equation]{Proposition}
\newtheorem{cor}[equation]{Corollary}
\newtheorem{lem}[equation]{Lemma}

\theoremstyle{definition}

\newtheorem{DEF}[equation]{Definition}
\newtheorem{rem}[equation]{Remark}

\newtheorem{tab}[equation]{Table}

\def\0b{\bar{0}}

\def\mod{\hbox{mod}}
\def\andd{\quad\hbox{and}\quad}
\def\ind{\hbox{ind}}

\def\v{{\mathcal V}}

\def\vd{\dot{\mathcal V}}
\def\vz{{\mathcal V}^{0}}
\def\vt{\tilde{\mathcal V}}

\def\fm{(\cdot,\cdot)}
\def\a{\alpha}

\def\ac{\alpha^{\vee}}
\def\bc{\betta^{\vee}}

\def\w{{\mathcal W}}

\def\sub{\subseteq}
\def\rd{\dot{R}}

\def\lam{\lambda}
\def\Lam{\Lambda}

\def\1k{\frac{1}{k}}

\def\la{\langle}
\def\ra{\rangle}

\def\d{\delta}

\def\b{\beta}
\def\bc{\beta^{\vee}}

\def\sg{\sigma}

\def\hh{{\mathcal H}}

\def\sgt{\tilde{\sg}}

\def\ti{{\mathcal T}_{i}}

\def\bbbz{{\mathbb Z}}
\def\bbbr{{\mathbb R}}

\def\da{\dot{\a}}

\def\proof{\noindent{\bf Proof. }}

\def\Lamt{\tilde{\Lam}}

\def\rank{\hbox{rank}}

\def\rr{\mathcal R}

\def\a{\alpha}

\def\ac{\alpha^{\vee}}

\def\andd{\quad\hbox{and}\quad}

\def\b{\beta}

\def\bc{\beta^{\vee}}

\def\d{\delta}

\def\e{\epsilon}

\def\hh{{\mathcal H}}

\def\bbbr{{\mathbb R}}

\def\ind{\hbox{ind}}
\def\fm{(\cdot,\cdot)}

\def\mod{\hbox{mod}}

\def\R{{\mathcal R}}
\def\sub{\subseteq}
\def\rd{\dot{R}}

\def\lam{\lambda}
\def\Lam{\Lambda}

\def\1k{\frac{1}{k}}

\def\la{\langle}
\def\ra{\rangle}

\def\GL{GL}

\def\qed{\hfill$\Box$\vspace{5mm}}
\def\sg{\sigma}

\def\rtimes{R^{\times}}

\def\sgt{\tilde{\sg}}

\def\ti{{\mathcal T}_{i}}

\def\Lamt{\tilde{\Lam}}

\def\vd{\dot{\mathcal V}}

\def\vz{{\mathcal V}^{0}}

\def\vt{\tilde{\mathcal V}}

\def\w{{\mathcal W}}

\def\bbbz{{\mathbb Z}}

\def\s{S}

\def\e{E}
\def\r{R}

\def\vb{\overline{\v}}

\def\ti{\tilde{I}}

\def\Lamt{\tilde{\Lam}}
\def\sgt{\tilde{\sg}}
\def\supp{\hbox{supp}}
\def\wdot{\dot{\w}}

\def\FO{\hbox{FO}}
\def\O{\hbox{O}}
\def\emp{\emptyset}

\begin{document}

\setcounter{page}{1}

\title{Simply laced extended affine Weyl groups\\
(a finite presentation)}

\author{Saeid Azam, Valiollah Shahsanaei}
\address
{Department of Mathematics\\ University of Isfahan\\Isfahan, Iran,
P.O.Box 81745-163.} \email{azam@sci.ui.ac.ir}
\thanks{The authors would like to thank the Center of Excellence for
Mathematics, University of Isfahan.} 

\begin{abstract} Extended affine Weyl groups are the Weyl groups of
extended affine root systems. Finite presentations for extended
affine Weyl groups are known only for nullities $\leq 2$, where
for nullity $2$ there is only one known such presentation. We give
a finite presentation for the class of simply laced extended
affine Weyl groups. Our presentation is nullity free if rank $>1$
and for rank $1$ it is given for nullities $\leq 3$. The
generators and relations are given uniformly for all types, and
for a given nullity they can be read from the corresponding finite
Cartan matrix and the semilattice involved in the structure of the
root system.
\end{abstract}

\maketitle \markboth{A FINITE PRESENTATION}{S. AZAM, V.
SHAHSANAEI}

\medskip
\setcounter{section}{-1}
\section{Introduction}\label{Introduction}
{\it Extended affine Weyl groups}  are the Weyl groups of {\it
extended affine root systems}  which are higher nullity
generalizations of affine and finite root systems. In 1985, K.
Saito [S] introduced axiomatically the notion of an extended
affine root system, he was interested in possible applications to
the study of singularities.

Extended affine root systems also arise as the root systems of a
class of infinite dimensional Lie algebras called {\it extended
affine Lie algebras}. A systematic study of extended affine Lie
algebras and their root systems is given in [AABGP], in particular
a set of axioms, different from those given by Saito \cite{S}, is
extracted from algebras for the corresponding root systems. In
[A2], the relation between axioms of [S] and [AABGP] for extended
affine root system's is clarified.

Let $R$ be an extended affine root system of nullity $\nu$ (see
Definition \ref{ears}) and $\v$ be the real span of $R$ which by
definition is equipped with a positive semi-definite bilinear form
$I$. We consider $R$ as a subset of a hyperbolic extension $\vt$
of $\v$, where $\vt$ is equipped with a non-degenerate extension
$\ti$ of $I$ (see Section \ref{hyperbolic}). Then the extended
affine Weyl group $\w$ of $R$ is by definition the subgroup of the
orthogonal group of $\vt$ generated by reflections based on the
set of non-isotropic roots $\rtimes$ of $R$.

This work is the first output of a three steps project on the
presentations of extended affine Weyl groups and its application
to the study of extended affine Lie algebras. In the first step,
we study finite presentations for extended affine Weyl groups,
where in this work we restrict ourself to the simply laced cases.
In the second step the results of the current work will apply to
investigate the existence of the so called a {\it presentation by
conjugation} for the simply laced extended affine Weyl groups (see
\cite{Kr} and \cite{A3}). Finally, in the third step, we will
apply the results of the second step to investigate validity of
certain classical results for the class of simply laced extended
affine Lie algebras. There is only a little known about the
presentations of extended affine Weyl groups. In fact if $\nu>2$,
there is no known finite presentation for this class and for
$\nu=2$ there is only one known finite presentation called the
{\it generalized Coexter presentation} (see \cite{ST}).

We give a finite presentation, for simply laced extended affine
Weyl groups, which is nullity free if $\rank>1$ and for $\rank=1$
it is given for nullities $\leq 3$ (see Theorem
\ref{presen-4-sim}). Our presentation highly depends on the
classification of semilattices (see Definition \ref{deflattice}),
up to similarity, which appears in the structure of extended
affine root systems (see (\ref{AABGP})). Since for types $A_\ell$
($\ell\geq 2$), $D_\ell$, $E_\ell$ for arbitrary nullity and for
type $A_1$ for nullity $\leq 3$ this classification is known (see
[AABGP, Chapter II]), our presentation is explicit for the
mentioned types and nullities.

The paper is arranged as follows. In Section \ref{hyperbolic}, we
obtain several results regarding the structure of certain
reflection groups. The results of this section are similar but
more general than those in \cite{MS} and \cite{A4}, and are
applicable to a wide range of reflection groups including extended
affine Weyl groups.
In Section \ref{heisenberg-sim},  we introduce a notion of {\it
supporting class} for the semilattices involved in the structure
of extended affine root systems. This notion plays a crucial role
in our work.  Also we study an intrinsic subgroup $\mathcal H$ of
an extended affine Weyl group $\w$ which we call it a {\it
Heisenberg-like group}. The center of $\hh$ is fully analyzed in
terms of the supporting class of the root system. This a basic
achievement which distinguishes the results of this section from
those in \cite{A4} and \cite{MS} (See Corollary \ref{centers-1}).
In particular it, together with other results, provides a unique
expression of Weyl group elements in terms of the elements of a
finite Weyl group, certain well-known linear transformations
belonging to the corresponding Heisenberg-like group and central
elements (see Proposition \ref{w-form}). We encourage the reader
to compare this with its similar results in \cite{MS} and
\cite{A4}.
The main results are given in Sections \ref{presentation-sim} and
\ref{presentation-W-sim} where we give our explicit presentations
for the extended affine Weyl group $\w$ and the Heisenberg-like
group $\mathcal H$.

The presentation is obtained as follows. First by analyzing the
semilattice involved in the structure of $R$, we obtain a finite
presentation for $\mathcal H$. The generators and relations depend
on nullity, the supporting class of the involved semilattice and
the Cartan matrix of the corresponding finite type. Next using the
fact that $\w=\dot{\w}\ltimes\mathcal H$, where $\dot{\w}$ is a
finite Weyl group of the same type of $R$, we obtain our
presentation for $\w$. So this presentation consists of three
parts, a presentation for $\mathcal H$, the Coxeter presentation
for the finite Weyl group $\dot{\w}$ and the relations imposed by
the semidirect product (see Theorems \ref{presen-1-sim} and
\ref{presen-4-sim}).

For a systematic study of extended affine Lie algebras and their
root systems we refer the reader to [AABGP]. For the study of
extended affine Weyl group we refer the reader to [S], [MS],
[A1,2,3,4], [ST] and [T].

\vspace{5mm}
\section{\bf REFLECTIONS GROUPS}\label{hyperbolic}
 \setcounter{equation}{0}
Let $\v$ be a finite dimensional real vector space equipped with a
non-trivial symmetric bilinear form $I=\fm$ of nullity $\nu$.   An
element $\a$ of $\v$  is called {\it nonisotropic} ({\it
isotropic}) if $(\a,\a)\not=0$ ($(\a,\a)=0$). We denote the set of
nonisotropic elements of a subset $A$ with $A^\times$. If $\a$ is
non-isotropic, we set $\ac=2\a/(\a,\a)$.  Let $\vz$ be the radical
of the form and $\vd$ be a fixed complement of $\vz$ in $\v$ of
dimension $\ell$. Throughout this section we fix a basis
$\{\a_1,\ldots,\a_\ell\}$ of $\vd$, a basis
$\{\sg_1,\ldots,\sg_\nu\}$ of $\v^0$ and a basis
$\{\lam_1,\ldots,\lam_\nu\}$ of $(\v^0)^*$. We enlarge the space
$\v$ to a $\ell+2\nu$-dimensional vector space as follows. Set
\begin{equation}\label{hyp}
\vt:=\v\oplus (\vz)^*,
\end{equation}
where  $(\vz)^*$ is the dual space of $\vz$. Now we extend the
from $\fm$ on $\v$ to a non-degenerate form, denoted again by
$\fm$, on $\vt$ as follows:
\begin{equation}\label{hyp1}
\begin{array}{ll}
 \bullet& (.,.)_{|_{\v\times\v}}:=(.,.),\vspace{2mm}\\
\bullet& (\vd,(\vz)^*)=((\vz)^*,(\vz)^*):=0,\vspace{2mm}\\
\bullet& (\sg_r,\lam_s):=\delta_{r,s},\quad 1\leq r,s\leq\nu
\end{array}
\end{equation}
The pair $(\vt,I)$ is called a {\it hyperbolic extension} of
$(\v,I)$.

Let $\O(\vt,I)$ be the orthogonal subgroup of $\GL(\vt)$, with
respect to $I=\fm$. We also set $$\FO(\vt,I)=\{w\in \O(\vt,I)\mid
w(\d)=\d\hbox{ for all }\d\in\v^0\}. $$  For $\a\in\v^\times$, the
element $w_\a\in \hbox{FO}(\vt,I)$ defined by
$$w_\a(u)=u-(u,\ac)\a,\qquad (u\in\vt),$$ is called the {\it
reflection based on} $\a$. It is easy to check that
\begin{equation}\label{normal}
ww_{\a}w^{-1}=w_{w(\a)}\qquad (w\in \O(\vt,I)).
\end{equation}
For a subset $A$ of $\v$, the group
\begin{equation}\label{reflection}
\w_A=\la  w_\a\mid\a\in A^\times\ra,
\end{equation}
is called the {\it reflection group in }$\hbox{FO}(\vt,I)$ {\it
based on} $A$.

For $\a\in\v$ and $\sg\in\vz$, we define $T^{\sg}_{\a}\in\hbox{\e
nd}(\vt)$ by
\begin{equation}\label{xi}
T^{\sg}_{\a}(u):=u-(\sg,u)\a+(\a,u)\sg
-\frac{(\a,\a)}{2}(\sg,u)\sg\quad(u\in\vt).
\end{equation}
The basic properties of the linear maps $T_\a^\sg$ are listed in
the following lemma. The terms of the form $[x,y]$ appearing in
the lemma denotes the commutator $x^{-1}y^{-1}xy$ of two elements
$x,y$ in a group.
\begin{lem}\label{new-4}
Let $\a,\b,\gamma\in\v$, $\sg,\d,\tau\in\v^0$ and
$w\in\hbox{FO}(\vt,I)$. Then

 (i) $T^{r\sg}_{\a}=T^{\sg}_{r\a}$, $r\in\bbbr$,

 (ii)
$T_{\a}^{\sg+\delta}=T_{\a}^{\sg}T_{\a}^{\delta}T_{\delta}^{(\a,\a)\frac{\sg}{2}}$,

(iii) $T^\sg_{\a+\b}=T^\sg_\a T^{\sg}_\b$,

(iv) $[T^{\sg}_\a,T^{\d}_\b]=T_\sg^{(\a,\b)\d}$,

(v)
$[T_{\a}^{\sg}T_{\b}^{\delta},T_{\gamma}^{\tau}]=[T_{\a}^{\sg},T_{\gamma}^{\tau}]
[T_{\b}^{\delta},T_{\gamma}^{\mu}]$,

(vi) $wT_{\a}^{\sg}w^{-1}=T_{w(\a)}^{\sg}$,

(vii) $ T^\sg_{\ac}=w_{\a+\sg}w_{\a}$, $\a\in\v^{\times},$

(viii) $(T^\d_\sg)^{-1}=T^\sg_\d$.
\end{lem}

\proof The proof of each statement can be seen using the
definition of the maps $T^\sg_\a$ and straightforward
computations.\qed

Part (vii) of Lemma \ref{new-4} has been in fact our motivation
for defining the maps $T_\a^\sg$. If we denote by $Z(G)$ the
center of a group $G$, then we have from parts (iii), (vi) and
(vii) of Lemma \ref{new-4} that for $\a\in\v$, $\sg,\d\in\v^0$,
\begin{equation}\label{neww}
T_\a^\sg\in\hbox{FO}(\vt,I)\andd T^\sg_\d\in
Z\big(\hbox{FO}(\vt,I)\big).
\end{equation}
\begin{lem}\label{formul-T}
Let $\a\in\v$, $n_r\in\bbbr$ and $1\leq r\leq\nu$. Then
$$T_{\a}^{\Sigma_{r=1}^{\nu}n_{r}\sg_{r}}=\prod_{r=1}^{\nu}T_{\a}^{n_{r}\sg_{r}}
\prod_{1\leq r<s\leq
\nu}T^{n_{r}n_{s}\frac{(\a,\a)}{2}\sg_{r}}_{\sg_s}.$$
\end{lem}

\proof For $ 1\leq r\leq\nu$, set
$\delta_{r}:=\sum_{i=r}^{\nu}n_{i}\sg_{i}$. Using Lemma
\ref{new-4}~(ii)-(iii) and (\ref{neww}), we have
\begin{eqnarray*}
\qquad\qquad\qquad
T_{\a}^{\Sigma_{r=1}^{\nu}n_{r}\sg_{r}}&=&T_{\a}^{n_{1}\sg_{1}+\delta_{2}}
=T_{\a}^{n_{1}\sg_{1}}T_{\a}^{\delta_{2}}
T_{\delta_{2}}^{n_1\frac{(\a,\a)}{2}\sg_{1}}\\
&=&T_{\a}^{n_{1}\sg_{1}}T_{\a}^{n_2\sg_{2}}T_{\a}^{\delta_{3}}
T_{\delta_{3}}^{n_{2}\frac{\sg_{2}}{2}}T_{\delta_{2}}^{n_{1}\frac{(\a,\a)}{2}\sg_{1}}\\
&\vdots&\\
&=&\prod_{r=1}^{\nu}T_{\a}^{n_{r}\sg_{r}}\prod_{r=1}^{\nu-1}
T_{\delta_{r+1}}^{n_{r}\frac{(\a,\a)}{2}\sg_{r}}\\
&=&\prod_{r=1}^{\nu}T_{\a}^{n_{r}\sg_{r}}\prod_{r=1}^{\nu-1}
\prod_{s=r+1}^{\nu}T_{\sg_{s}}^{n_rn_s\frac{(\a,\a)}{2}\sg_{r}}\\
\qquad\qquad\qquad&=&\prod_{r=1}^{\nu}T_{\a}^{n_{r}\sg_{r}}\prod_{1\leq
r<s\leq
\nu}T^{n_rn_s\frac{(\a,\a)}{2}\sg_{r}}_{\sg_s}.\hbox{\qquad\qquad\qquad\qed}
\end{eqnarray*}

For a subset $A$ of $\v$, consider the subgroup
\begin{equation}\label{hb}
{\mathcal H}(A):=\la  w_{\a+\sg}w_\a\mid \a\in A^\times,
\;\sg\in\vz,\;\a+\sg\in A\ra.
\end{equation}
of $\w_A$. We note that ${\mathcal H}(A)\leqslant {\mathcal
H}(\v)=\la w_{\a+\sg}w_\a\mid \a\in \v^\times, \;\sg\in\vz\ra$.

 We recall
that a group $H$ is called {\it two-step nilpotent} if the
commutator $[H,H]$ is contained in the center $Z(H)$ of $H$. We
also recall that in a two-step nilpotent group $H$, the commutator
is bi-multiplicative, that is
\begin{equation}\label{bicomm}
[\prod_{i=1}^{n}x_{i},\prod_{j=1}^{m}y_{j}]=\prod_{i=1}^{n}\prod_{j=1}^{m}
  [x_{i},y_{j}]
  \end{equation}
for all $x_{i},y_{j}\in H.$
\begin{lem}\label{xi-4}
Let  $A$ be a subset of $\v$. Then

 (i) ${\mathcal H}(A)$ is a two-step
nilpotent group.

(ii) If $X$ is a generating subset of ${\mathcal H}(A)$, then
$[X,X]$ generates the commutator $[{\mathcal H}(A),{\mathcal
H}(A)]$ of ${\mathcal H}(A)$.
\end{lem}
\proof (i) If $w_\a w_{\a+\sg}$ and $w_\b w_{\b+\d}$ are two
generators of ${\mathcal H}(A)$, $\a,\b\in A^\times$,
$\a+\sg,\b+\d\in A$, then by Lemma \ref{new-4}(vii),(iv), we have
$$
[w_{\a+\sg}w_\a,w_{\b+\d}w_\b]=[T_{\ac}^\sg,T_{\bc}^{\d}]=T_\sg^{(\ac,\bc)\d}.
$$ Now the result follows from (\ref{neww}). (ii) is an immediate
consequence of (i) and (\ref{bicomm}).\qed

\begin{lem}\label{xi-5} ${\mathcal H}(\v)=\la
T_{\a_i}^{n_{i,r}\sg_r}, T_{\sg_r}^{m_{r,s}\sg_s}\mid 1\leq
i\leq\ell;\; 1\leq r,s\leq\nu;\; n_{i,r}, m_{r,s}\in\bbbr\ra.$
\end{lem}

\proof From  parts (iii) and (vii) of Lemma \ref{new-4} and Lemma
\ref{formul-T} we see that ${\mathcal H}(\v)$ is a subset of the
right hand side. Conversely, it follows from Lemma
\ref{new-4}(vii) and (iv) that the right hand side is a subset of
${\mathcal H}(\v)$. \qed

\begin{lem}\label{new-xi-6}
Let
$h=\prod_{r=1}^{\nu}\prod_{i=1}^{\ell}T_{\a_i}^{n_{i,r}\sg_r}\prod_{1\leq
r<s\leq\nu} T_{\sg_r}^{m_{r,s}\sg_s}$, where $n_{i,r},
m_{r,s}\in\bbbr$. If $\b_j=\sum_{i=1}^{\ell}n_{i,j}\a_i$, $1\leq
j\leq\nu$, then $$\begin{tabular}{c}
   $h(\lam_j)=\lam_r- \b_j-
\frac{(\b_j,\b_j)}{2}\sg_j-\sum_{ 1\leq r\leq
j-1}m_{r,j}\sg_r+\sum_{ j+1\leq s\leq\nu}m_{j,s}\sg_s.$
 \end{tabular}$$
\end{lem}

\proof Let $1\leq  j\leq\nu$ and $\da\in\vd$. Then  from
(\ref{xi}) and (\ref{hyp1}), it follows that
 $$T_{\da}^{\sg_j}(\lam_j)=\lam_j-
\da-\frac{(\da,\da)}{2}\sg_j\andd
T_{\sg_r}^{\sg_s}(\lam_j)=\lam_j-\delta_{s,j}\sg_r+\delta_{r,j}\sg_s,\quad
1\leq r,s\leq\nu,$$
 and so using  (\ref{neww}) and
Lemma \ref{new-4}(ii)-(iii), we have
\begin{eqnarray*}
h(\lam_j)&=&\prod_{r=1}^{\nu}\prod_{i=1}^{\ell}T_{\a_i}^{n_{i,r}\sg_r}\prod_{1\leq
r<s\leq\nu} T_{\sg_r}^{m_{r,s}\sg_s}(\lam_j)\\
&=&\prod_{r=1}^{\nu}
\prod_{i=1}^{\ell}T_{\a_i}^{n_{i,r}\sg_r}\prod_{j+1\leq s\leq\nu}
T_{\sg_j}^{m_{j,s}\sg_s}\prod_{ 1\leq r\leq j-1}
T_{\sg_r}^{m_{r,j}\sg_j}\hspace{-5mm}\prod_{\{1\leq r<s\leq\nu\mid
r,s\neq j\}} \hspace{-5mm}T_{\sg_r}^{m_{r,s}\sg_s}(\lam_j)\\
&=&\prod_{r=1}^{\nu} \prod_{i=1}^{\ell}T_{\a_i}^{n_{i,r}\sg_r}
T_{\sg_j}^{\sum_{j+1\leq s\leq\nu}m_{j,s}\sg_s}
T_{\sum_{ 1\leq r\leq j-1}m_{r,j}\sg_r}^{\sg_j}(\lam_j)\\
&=&\prod_{r=1}^{\nu}
\prod_{i=1}^{\ell}T_{\a_i}^{n_{i,r}\sg_r}T_{\sg_j}^{\sum_{j+1\leq
s\leq\nu}m_{j,s}\sg_s}(\lam_j-\sum_{ 1\leq r\leq
j-1}m_{r,j}\sg_r)\\
&=&\prod_{r=1}^{\nu}
\prod_{i=1}^{\ell}T_{\a_i}^{n_{i,r}\sg_r}(\lam_j- \sum_{ 1\leq
r\leq j-1}m_{r,j}\sg_r+\sum_{j+1\leq s\leq\nu}m_{j,s}\sg_s)\\
&=&\prod_{r=1}^{\nu}
\prod_{i=1}^{\ell}T_{\a_i}^{n_{i,r}\sg_r}(\lam_j)- \sum_{ 1\leq
r\leq j-1}m_{r,j}\sg_r+\sum_{j+1\leq s\leq\nu}m_{j,s}\sg_s\\
\end{eqnarray*}
\begin{eqnarray*}
 &=&\prod_{i=1}^{\ell}T_{\a_i}^{n_{i,j}\sg_j}(\lam_j)- \sum_{ 1\leq
r\leq j-1}m_{r,j}\sg_r+\sum_{j+1\leq s\leq\nu}m_{j,s}\sg_s\\
&=&T_{\b_j}^{\sg_j}(\lam_j)-\sum_{ 1\leq r\leq j-1}m_{r,j}\sg_r+
\sum_{j+1\leq s\leq\nu}m_{j,s}\sg_s\\ \qquad\qquad&=&\lam_j- \b_j-
\frac{(\b_j,\b_j)}{2}\sg_j-\sum_{ 1\leq r\leq j-1}m_{r,j}\sg_r
+\sum_{j+1\leq s\leq\nu}m_{j,s}\sg_s.\hbox{\qquad\qquad\qed}
\end{eqnarray*}
\begin{lem}\label{xi-6}
Each element $h\in {\mathcal H}(\v)$ has a unique expression in
the form
\begin{equation}\label{exp} h=h(n_{i,r},m_{r,s}):=
\prod_{r=1}^{\nu}\prod_{i=1}^{\ell}T_{\a_i}^{n_{i,r}\sg_r}\prod_{1\leq
r<s\leq\nu} T_{\sg_r}^{m_{r,s}\sg_s}\quad(n_{i,r},
m_{r,s}\in\bbbr).
\end{equation}
\end{lem}
\proof  Let $h\in {\mathcal H}$. From (\ref{neww}), Lemmas
\ref{xi-5} and \ref{new-4}(iv) it follows that $h$ has an
expression in the form (\ref{exp}). Now let $h(n'_{i,r},m'_{r,s})$
be another expression of $h$ in the form (\ref{exp}). Then by
acting these two expressions of $h$ on $\lam_j$'s, $1\leq
j\leq\nu$, we get from Lemma \ref{new-xi-6} that
$n_{i,r}=n'_{i,r}$ and $m_{r,s}=m'_{r,s}$ for all  $1\leq i\leq
\ell$ and $1\leq r<s\leq\nu$.\qed

\begin{lem}\label{xi-7}
(i) $Z\big({\mathcal H}(\v)\big)=\la T_{\sg_r}^{m_{r,s}\sg_s}\mid
1\leq r<s\leq\nu,\;m_{r,s}\in\bbbr\ra.$

(ii) For any fixed nonzero real numbers $m_{r,s}$, $1\leq
r<s\leq\nu$, the group $\la T_{\sg_r}^{m_{r,s}\sg_s}\mid 1\leq
r<s\leq\nu\ra$ is  free abelian of rank $\frac{\nu(\nu-1)}{2}$ .

(iii) ${\mathcal H}(\v)$ is a torsion free group.
\end{lem}

\proof (i) By (\ref{neww}), Lemmas \ref{new-4}(vii)  and
\ref{xi-5}, it is clear that the right hand side in the statement
 is a subset of
the left hand side. To show the reverse inclusion, let $h\in
Z({\mathcal H}(\v))$. Consider an expression $h(n_{i,r},m_{r,s})$
of $h$ in the form (\ref{exp}). We must show that $n_{i,r}=0$, for
all $1\leq i\leq\ell$ and $1\leq r\leq\nu$. Since ${\mathcal
H}(\v)$ is a two-step nilpotent group, we have from (\ref{bicomm})
and Lemmas \ref{xi-5} and \ref{new-4}(iv) that for all $1\leq
j\leq\ell$ and $1\leq s\leq \nu$,
\begin{eqnarray*}
1=[h,T_{\a_j}^{\sg_s}]&=&\prod_{r=1}^{\nu}\prod_{i=1}^{\ell}
 [T_{\a_i}^{n_{i,r}\sg_r},T_{\a_j}^{\sg_s}]\\
 &=&
 \prod_{r=1}^{\nu}\prod_{i=1}^{\ell}T_{\sg_r}^{n_{i,r}(\a_{i},\a_{j})\sg_s}=
 \prod_{r=1}^{\nu}T_{\sg_r}^{(\sum_{i=1}^{\ell}n_{i,r}\a_{i},\a_{j})\sg_s}.
\end{eqnarray*}
Therefore by Lemma \ref{xi-6},
$\sum_{i=1}^{\ell}n_{i,r}(\a_{i},\a_{j})=0$, for all $1\leq
j\leq\ell$ and $1\leq r\leq\nu$. But
$\vd=\sum_{i=1}^{\ell}\bbbr\a_{i}$ and the form restricted to
$\vd$ is non-degenerate, so $n_{i,r}=0$ for all $1\leq i\leq\ell$,
$1\leq r\leq\nu$.

(ii) We show that $\{T^{m_{r,s}\sg_s}_{\sg_r}\mid 1\leq r\leq
s\leq\nu\}$ is a free basis for the group under consideration. Let
$\prod_{1\leq r<s\leq\nu}^{\nu}T^{n_{r,s}m_{r,s}\sg_s}_{\sg_r}=1$,
$n_{r,s}\in\bbbz$ for $1\leq r<s\leq\nu$. Then by Lemma
\ref{xi-6}, $n_{r,s}m_{r,s}=0$ for all $r,s$. The result now
follows as $m_{r,s}$'s are nonzero.

(iii) Let $h\in {\mathcal H}(\v)$ and assume $h^n=1$ for some
$n\in\mathbb N$. Let $h=h(n_{i,r},m_{r,s})$ (see (\ref{exp})).
Since ${\mathcal H}(\v)$ is two-step nilpotent, we have
\begin{eqnarray*}
1=h^n=\prod_{r=1}^{\nu}\prod_{i=1}^{\ell}T_{\a_i}^{nn_{i,r}\sg_r}c\prod_{1\leq
r<s\leq\nu} T_{\sg_r}^{nm_{r,s}\sg_s},
\end{eqnarray*}
where $c$ is a central element. By part (i),
$1=h^n=h^n(nn_{i,r},m'_{r,s})$ for some real numbers $m'_{r,s}$.
By Lemma \ref{xi-6}, $n_{i,r}=0$ for all $i,r$. Thus $h$ is
central and so  $c=1$. Now it follows again from Lemma \ref{xi-6}
that $m_{r,s}=0$ for all $r,s$.\qed

\begin{cor}\label{xi-9} For any subset $A$ of $\v$, ${\mathcal H}(A)$
is a torsion free group.
\end{cor}

 Recall that we
have fixed a complement $\vd$ of $\v^0$ in $\v$. Now for a subset
$A$ of $\v$, we set $$\dot{A}:=\{\a\in\vd\mid \a+\sg\in A\hbox{
for some
 }\sg\in\v^0\}.$$

 \begin{pro}\label{xi-10}
 Let $A$ be a subset of $\v$ such that $w_\a(A)\sub A$ for all $\a\in A^\times$ and
 ${\mathcal H}(A)=\la
w_{\a+\sg}w_\a\mid\a\in \dot{A}^\times,\sg\in\v^0,\a+\sg\in A\ra$.
Then $\w_A=\w_{\dot{A}}\ltimes {\mathcal H}(A).$
\end{pro}
\proof  Let $\a\in{A}$, $\sg\in\v^0$, ${\a}+\sg\in A$ and
$w\in\w_A$. By assumption  $w(\a)$ and $w(\a+\sg)=w(\a)+\sg$ are
elements of $A$. Thus, ${\mathcal H}(A)$ is a normal subgroup of
$\w_A$, by (\ref{normal}). Now we show that ${\mathcal H}(A)\cap
\w_{A}=\{1\}$. Since ${\mathcal H}(A)\sub {\mathcal H}(\v)$ and
$\w_{\dot{A}}\subseteq\w_{\vd}$, it is enough to show that
${\mathcal H}(\vd)\cap \w_{\vd}=\{1\}$. Let $h\in\w_{\dot{\v}}\cap
{\mathcal H}(\v)$ and consider an expression of $h$ in the form
(\ref{exp}). Since $h\in\w_{\dot{\v}}$, we have from (\ref{hyp1})
that $h(\lam_j)=\lam_j$, $1\leq j\leq\nu$ and so it follows from
Lemma \ref{new-xi-6} that $h=1$.

To complete the proof, we must show $\w_A=\w_{ \dot{A}} {\mathcal
H}(A)$. Let $\a\in \dot{A}$ and $\sg\in\vz$ such that $\da+\sg\in
A$. By assumption, $w_{\da+\sg}w_{\da}\in {\mathcal
H}(A)\leqslant\w_A$ and $w_{\da+\sg}\in\w_A$, therefore
$w_{\da}=w_{\da+\sg}w_{\da+\sg}w_{\da}\in\w_A$ and so $\w_{
\dot{A}}\sub\w_A$. This shows that $\w_{\dot{A}} {\mathcal
H}(A)\sub\w_A$. To see the reverse inclusion, let $w_{\a}$, $\a\in
A^\times$, be a generator of $\w_A$. Then $\a=\da+\sg$ where
$\da\in \dot{A}^{\times}$ and $\sg\in\vz$. Since
$w_{\da}\in\w_{\dot{A}}$ and since by assumption
$w_{\da}w_{\da+\sg}\in {\mathcal H}(A)$, we have
$w_\a=w_{\da+\sg}=w_{\da}(w_{\da}w_{\da +\sg}) \in
\w_{\dot{A}}{\mathcal H}(A)$. This completes  the proof.\qed

\vspace{5mm}
\section{\bf EXTENDED AFFINE WEYL GROUPS }\label{heisenberg-sim}
 \setcounter{equation}{0}

In this section, we study Weyl groups of simply laced extended
affine root systems. We are mostly interested in finding a
particular finite
 set of generators for such a Weyl group and its center (see Proposition
 \ref{propo-impor}). Since the semilattice involved in the
 structure of an extended affine root system plays a crucial role
 in our study, we start this section with recalling the definition of
 a semilattice from \cite[II.\S 1]{AABGP} and introducing
 a notion of {\it supporting class}
for semilattices.  For the theory of extended affine root systems
the reader is referred to [AABGP]. In particular, we will use the
notation and concepts introduced there without further
explanations.

\begin{DEF}\label{deflattice}
A {\it semi-lattice} is a subset $\s$ of a finite dimensional real
vector space $\vz$ such that $0\in S$, $\s\pm2\s\subseteq \s$, $S$
spans $\vz$ and $S$ is discrete in $\vz$. The {\it rank }of $S$ is
defined to be the dimension $\nu$ of $\vz$.  Note that the
replacement of $\s\pm2\s\subseteq \s$ by $\s\pm\s\subseteq \s$ in
the definition gives one of the equivalent definitions for a {\it
lattice} in $\vz$. Semilattices $S$ and $S'$  in $\vz$ are said to
be similar if there exist $\psi\in GL(\vz)$ so that
$\psi(S)=S'+\sg'$ for some $\sg'\in S'$.
  \end{DEF}
Let $S$ be a semilattice in $\v^0$ of rank $\nu$.
The $\bbbz$-span $\Lam$ of $S$ is a lattice in $\v^0$, a free
abelian group of rank $\nu$ which has an $\bbbr$-basis of $\v^0$
as its $\bbbz$-basis. By \cite[ II.1.11] {AABGP}, $S$ contains a
subset $B=\{\sg_1,\ldots,\sg_\nu\}$ of $S$ which forms a basis for
$\Lam$. We call such a set $B$, a {\it basis} for $S$. Then
$$\Lam=\la S\ra=\sum_{i=1}^\nu\bbbz\sg_r\;\;\hbox{ with }
\;\;\sg_r\in S\;\;\hbox{ for all }\;\;r.
$$
Consider $\Lamt:=\Lam/2\Lam$ as a $\bbbz_2$-vector
space with ordered basis $\tilde{B}$, the image of $B$ in $\Lamt$.
For $\sg\in\Lam$ and $1\leq r\leq\nu$, let $\sg(r)\in\{0,1\}$ be
the unique integer such that $\sgt=\sum_{r=1}^{\nu}\sg(r)\sgt_r$.
Then we set
$$\mbox{supp}_B(\sg):=\{1\leq r\leq\nu\mid \sg(r)=1\}.$$ Then
$\sg=\sum_{r\in\supp_{B}(\sg)}\sg_r$ ($\mod\; 2\Lam$). By [AABGP,
II.1.6], $S$ can be written in the form
$$
S=\bigcup_{j=0}^{m}(\d_j+2\Lam), $$ where $\d_0=0$ and $\d_j$'s
are distinct coset representatives of $2\Lam$ in $\Lam$. The
integer $m$ is called the {\it index}  of $S$ and is denoted by
$\ind(S)$.  The collection
\begin{equation}\label{support}
\supp_B(S):=\big\{\supp_B(\d_j)\mid 0\leq j\leq m\big \}
\end{equation}
is called the {\it supporting class of} $S$, with respect to $B$.
Since $\d_j=\sum_{r\in \supp_B(\d_j)}\sg_r$ ($\mod\; 2\Lam$), the
supporting set determines $S$ uniquely. Therefore, we may write
\begin{equation}\label{unique-1}
S=\biguplus_{J\in\supp_B(S)}(\tau_{_J}+2\Lam)\quad\hbox{where}\quad
\tau_{_J}:=\sum_{r\in J}\sg_r.
\end{equation}
 (By convention we have
$\tau_{_{\emptyset}}:=\sum_{r\in \emptyset}\sg_r=0$).  By [A5,
Proposition 1.12], if $\nu\leq 3$, then the index determines
uniquely, up to similarity, the semilattices in $\Lam$. So by
\cite[Table II.4.5]{AABGP}, up to similarity, the semilattices of
rank $\leq 3$ in $\Lam$ are listed in the following table,
according to their supporting classes:

\pagebreak
\begin{tab}\label{tab-1}
The supporting classes of semilattices, up to similarity, for
$\nu\leq 3$.
 \end{tab}
 $\begin{tabular}{c|c|c}
    $\nu$ &$\hbox{index}$ & $\supp_B(S)$ \\
    \hline
  0 & 0 & $\{\emp\}$ \\
   \hline
  1 & 1 & $\{\emp,\{1\}\}$ \\
  \hline
  2 & 2 &  $\{\emp,\{1\},\{2\}\}$\\
   & 3&  $\{\emp,\{1\},\{2\},\{1,2\}\}$\\
   \hline
  3 & 3 & $\{\emp,\{1\},\{2\},\{3\}\}$\\
    & 4 & $\{\emp,\{1\},\{2\},\{3\},\{2,3\}\}$\\
    & 5 & $\{\emp,\{1\},\{2\},\{3\},\{1,3\},\{2,3\}\} $\\
       & 6 & $\{\emp,\{1\},\{2\},\{3\},\{1,2\},\{1,3\},\{2,3\}\}$\\
      & 7&$\{\emp,\{1\},\{2\},\{3\},\{1,2\},\{1,3\},\{2,3\},\{1,2,3\}\}$\\
    \hline
\end{tabular}$
\vspace{.5cm}

Next we recall the definition of an extended affine root system.

\begin{DEF}\label{ears}
{\rm A subset $R$ of a non-trivial finite dimensional real vector
space $\v$, equipped with a positive semi-definite symmetric
bilinear form $(\cdot,\cdot)$, is called an {\it extended affine
root system} if $\r$ satisfies the following 8 axioms:
\begin{itemize}

\item  R1) $ 0\in \r$,

\item R2) $ -\r=\r$,

\item R3) $\r$ spans $\v$,

\item R4) $\a \in \r^{\times} \Longrightarrow 2\a \notin \r$,

\item R5) $\r$ is discrete in $\v$,

\item R6) For $\a \in \r^{\times}$  and $\beta\in \r$, there exist
non-negative integers $d,u$ such that $\b+n\a\in R$, $n\in\bbbz$,
if and only if $-d\leq n\leq u$, moreover $(\b,\ac)=d-u$,

\item R7) If $\r=\r_{1}\cup \r_{2}$, where $(\r_{1},\r_{2})=0$,
then either $\r_{1}=\emptyset$ or $\r_{2}=\emptyset$,

\item  R8)   For any $\sg \in \r^{0}$, there exists $\a\in
\r^{\times}$ such that $\a+\sg\in R$.
\end{itemize}

The dimension $\nu$ of the radical $\v^0$ of the form is called
the {\it nullity} of $R$, and the dimension $\ell$ of
$\vb:=\v/\v^0$ is called the {\it rank} of $R$. Sometimes, we call
$R$ a $\nu$-extended affine root system.}   Corresponding to the
integers $\ell$ and $\nu$, we set
$$J_\ell=\{1,\ldots,\ell\}\andd J_\nu=\{1,\ldots,\nu\}.
$$
\end{DEF}

Let $\r$ be a $\nu$-extended affine root system. It follows that
the form restricted to $\vb$ is positive definite and that
$\bar{\r}$, the image of $R$ in $\vb$, is an irreducible finite
root system (including zero) in $\vb$ ([AABGP, II.2.9]). The {\it
type} of $R$ is defined to be the type of the finite root system
$\bar{R}$. In this work {\it we always assume that $R$ is an
extended affine root system of simply laced type}, that is it has
one of the types $X_\ell=A_\ell$, $D_\ell,$ $E_6,$ $E_7$ or $E_8.$

According to [AABGP, II.2.37], we may fix a complement $\vd$ of
$\v^0$ in $\v$ such that
\begin{equation}\label{rd}
\rd:=\{\da\in\vd\mid\da+\sg\in R\hbox{ for some }\sg\in\v^0\}
\end{equation}
is a finite root system in $\vd$, isometrically isomorphic to
$\bar{R}$, and that
 $R$ is of the form
\begin{equation}\label{AABGP}
R=R(X_\ell,S)=(S+S)\cup(\rd+S)
\end{equation} where $S$ is
a semilattice in $\v^0$. Here $X_\ell$ denotes the type of $\rd$.
It is known that if $\ell>1$, then $S$ is a lattice in $\vz$.

Throughout our work, we fix two sets
$$
\Pi=\{\a_1,\ldots,\a_{\ell}\}\andd B=\{\sg_1,\ldots,\sg_\nu\},
$$
where $\Pi$ is a fundamental system for $\rd$ with $(\a_i,\a_i)=2$
for all $i\in J_\ell$, and $B$ is a basis for $S$. In particular,
$S$ has the expression as in (\ref{unique-1}). Since $B$ is fixed,
we write $\supp(S)$ instead of $\supp_B(S)$. From $S\pm
2S\subseteq S$, it follows that $\bbbz\sg_r\sub S$ for $r\in
J_\nu$, and so we have from (\ref{AABGP}) that
\begin{equation}\label{roots-sim}
\a_i+\bbbz\sg_r\subseteq \r,\qquad (i,r)\in J_\ell\times J_\nu
\end{equation}

 As in Section \ref{hyperbolic}, let
$\vt=\vd\oplus\v^0\oplus(\v^0)^*$, where $\vd$ is the real span of
$\rd$. With respect to the basis $B$, we extend the from $\fm$ on
$\v$ to a non-degenerate form, denoted again by $\fm$, on $\vt$ by
(\ref{hyp1}).

 We recall from (\ref{reflection}) and (\ref{hb}) that
 \begin{equation}\label{reform0-sim}
\w_R=\la w_{\a}\mid\a\in\rtimes\ra
\end{equation}
is a subgroup of $\FO(\vt,I)$ and
$$\mathcal H(R)=\la
w_{\a+\sg}w_\a\mid\a\in\rtimes,\sg\in\v^0,\a+\sg\in R\ra$$ is a
subgroup of $W_R$.

\begin{DEF}\label{EAWG} The  groups $\w_R$ and  $\mathcal H(R)$
  are called the
{\it extended affine Weyl group} and the {\it Heisenberg-like
group} of $R$, respectively. Since $\rd\sub R$, we may identify
the finite Weyl group of $\rd$ with the subgroup $\wdot=\la
w_\a\mid\a\in\rd^\times\ra$ of $\w_R$. When there is no confusion
we simply write $\w$ and ${\mathcal H}$ instead of $\w_R$ and
$\mathcal H(R)$ respectively.
\end{DEF}
\begin{lem}\label{reform2-sim}  ${\mathcal
H}=\la T^{\sg}_{\da}\mid\da\in\rd,\;\sg\in S\ra.$
\end{lem}

\proof For $\da\in\rd$ and $\sg\in S$, we have from Lemma
\ref{new-4}(vii) and (\ref{AABGP}) that
$T^{\sg}_{\da}=w_{\da+\sg}w_{\da}\in\mathcal H$. Also if
$\a\in\rtimes$, $\sg\in\v^0$ and $\a+\sg\in R$, then
$\a=\da+\tau$, where $\da\in\rd$ and by (\ref{AABGP}),
$\tau,\sg,\tau+\sg\in S$. Then
$w_{\a+\sg}w_\a=w_{\da+\tau+\sg}w_{\da}w_{\da}w_{\da+\tau}=
T_{\da}^{\tau+\sg}T_{-\da}^{\tau}.$ \qed

We next want to find certain finite sets
 of generators for both $\w$ and $\mathcal H$ and  their centers.
For $r,s\in J_\nu$, we set
\begin{equation}\label{def-crs}
c_{r,s}:=T_{\sg_{r}}^{\sg_{s}}\andd C:=\la c_{r,s}\mid 1\leq
r<s\leq\nu\ra.
\end{equation}
Then by (\ref{neww}) and  Lemma \ref{new-4}(viii) for all $r,s\in
J_\nu$, we have
\begin{equation}\label{comute}
c_{r,r}=c_{s,r}c_{r,s}=1\andd C\leq Z\big(\FO(\vt,I)\big).
\end{equation}
Moreover from  (\ref{def-crs}) and Lemma \ref{xi-7}(ii), it
follows that
\begin{equation}\label{free-abelian}
\mbox{$C$ is a free abelian group of rank $\nu(\nu-1)/2.$}
\end{equation}
Also for $(i,r)\in J_\ell\times  J_\nu$, we set
\begin{equation}\label{def-tir}
t_{i,r}:=T_{\a_{i}}^ {\sg_{r}}.
\end{equation}
From  (\ref{roots-sim}) and parts (iv) and (vii) of Lemma
\ref{new-4} for all $r,s\in J_\nu$ and $i,j\in J_\ell$,  one can
see easily that
\begin{equation}\label{comu-sim}
t_{i,r}=w_{\a_i+\sg_r}w_{\a_i}\in {\mathcal H}\andd
[t_{i,r},t_{j,s}]=c_{r,s}^{(\a_{i},\a_{j})}\in {\mathcal H}.
\end{equation}
\begin{lem}\label{wtw-sim}
  $w_{\a_{i}}t_{j,r}w_{\a_{i}}=
t_{j,r}t_{i,r}^{-(\a_i,\a_j)}$, for $i,j\in J_\ell$, $r\in J_\nu$.
 \end{lem}
 \proof Let $i,j\in J_\ell$, $r\in J_\nu$ and $\a=w_{\a_i}(\a_j)=\a_j-(\a_j,\a_i)\a_i$.
  We have from
  Lemma \ref{new-4}(vi) that
 $ w_{\a_{i}}t_{j,r}w_{\a_{i}}=T_{\a_i}^{\sg_{r}}=
 t_{j,r} t_{i,r}^{-(\a_i,\a_j)}.$
 \qed

In order to describe the centers $Z(\mathcal H )$ and $Z(\w)$, we
use the notion of supporting class of $S$ (with respect to $B$) by
assigning a subgroup of $C$ to $S$ as follows. We set
\begin{equation}\label{def-F(S)}
F(S):=\la z_{_J}\mid J\sub J_\nu\ra\leq C,
\end{equation}
where
 \begin{equation}\label{def-zJ}
z_{_J}:=\left\{\begin{array}{ll} \prod_{\{r,s\in J\mid r<s\}}
  c_{_{r,s}} & \hbox{if }J\in\supp(S)\vspace{3mm}\\
   c_{_{r,s}}^{2} &  \hbox{if
   $J=\{r,s\}\not\in\supp(S)$,}\vspace{3mm}\\
   1, & \hbox{otherwise}.
\end{array}\right.
\end{equation}
(Here we interpret the product on an empty index set to be $1$.)
Our goal is to prove that $Z(\w)=Z(\mathcal H )=F(S)$. Note that
if $\{r,s\}\sub\supp(S)$, then $z_{_{\{r,s\}}}=c_{r,s}$. In
particular, if $S$ is a lattice then the second condition in the
definition of $z_{_J}$ is surplus and so $F(S)=\la
z_{_{\{r,s\}}}\mid 1\leq r<s\leq\nu\ra=C$. Also from the way
$z_{_{\{r,s\}}}$ is defined and (\ref{comu-sim}) we note that
\begin{equation}\label{commutator-sim}
[t_{i,r},t_{j,s}]=c_{r,s}^{(\a_i,\a_j)}\in \la
z_{_{\{r,s\}}}\ra,\quad i,j\in J_\ell,\quad r,s\in J_\nu.
\end{equation}
\begin{lem}\label{nice-lemma}
Let $\a=\sum_{i=1}^{\ell}m_{i}\a_{i}\in\rd$ and
$\sg=\sum_{r=1}^{\nu}n_{r}\sg_{r}\in \Lam$. Then $$T_{\a}^{\sg}=
\prod_{r=1}^{\nu} \prod_{i=1}^{\ell} t_{i,r}^{m_in_r} \prod_{1\leq
r<s\leq\nu}c_{s, r}^{n_{r}n_{s}}.$$
\end{lem}

\proof  Using Lemma \ref{formul-T} and  Lemma \ref{new-4}(iii), we
have
\begin{eqnarray*}
T_{\a}^{\sg}&=&
 \prod_{r=1}^{\nu}T_{\a}^{n_{r}\sg_{r}}
\prod_{1\leq r<s\leq\nu}(T_{\sg_s}^{\sg_r})^{n_{s}n_{r}}=
\prod_{r=1}^{\nu} \prod_{i=1}^{\ell} t_{i,r}^{m_in_r} \prod_{1\leq
r<s\leq\nu}c_{s, r}^{n_{r}n_{s}}.
\end{eqnarray*}
\qed

For a subset  $J=\{i_1,\ldots, i_n\}$ of $J_\nu$ with $
i_1<i_2<\cdots <i_n$ and a group $G$ we make the convention
$$\prod_{i\in J} a_i=a_{i_1}a_{i_2}\cdots a_{i_n}\qquad\qquad(a_i\in G). $$
\begin{pro}\label{gen-H}
${\mathcal H}=\la t_{i,r},z_{_J}\mid (i,r)\in J_\ell\times
J_\nu,\;J\subseteq J_\nu\ra$.
\end{pro}

\proof Let $T$ be the group in the right hand side of the
equality. We proceed with the proof in the following two steps.

(1) $T\subseteq {\mathcal H}$. By (\ref{comu-sim}), it is enough
to show that $z_{_J}\in {\mathcal H}$ for any set $J\sub J_\nu$.
First, let $J\in\supp(S)$.  Then by the definition of $z_{_J}$, we
have $z_{_J}=\prod_{\{r,s\in J\mid r<s\}}c_{r,s}$.  Now it follows
from Lemma \ref{formul-T} that
  \begin{eqnarray*}
  T_{\a_{i}}^{\tau_{_J}}=
  T_{\a_{i}}^{\Sigma_{r\in J}\sg_{r}}
  &=&\prod_{\{r,s\in J\mid r<s\}}T^{\sg_r}_{\sg_s}
 \prod_{r\in J}T_{\a_{i}}^{\sg_{r}}=\prod_{\{r,s\in J\mid r<s\}}c_{s,r}
 \prod_{r\in J}T_{\a_{i}}^{\sg_{r}}=z_{_J}^{-1}\prod_{r\in
 J}t_{i,r}.
  \end{eqnarray*}
 But since $\a_{i}+\tau_{_J}\in\r$ (by (\ref{unique-1}) and (\ref{AABGP})), it follows from
Lemma \ref{reform2-sim} and (\ref{comu-sim}) that
$$z_{_J}=(T_{\a_{i}}^{\tau_{_J}})^{-1}\prod_{r\in
    J}t_{i,r}\in {\mathcal H}.$$
Finally, suppose $J=\{r,s\}\not\in\supp(S)$ where $1 \leq r
<s\leq\nu$. Then from the definition of $z_{_J}$ and
(\ref{comu-sim}) we have
$$z_{_J}=c_{r,s}^{2}=[t_{i,r},t_{i,s}]\in{\mathcal H}.$$
 This completes the proof of step (1).

 (2) ${\mathcal H}\subseteq  T$. We have from Lemmas
 \ref{reform2-sim}, \ref{new-4}(ii) and (\ref{unique-1}) that
\begin{eqnarray*} {\mathcal H}&=&
\la T_{\a}^{\sg}\mid \a\in\rd,\; \sg\in\s\ra\\&=&\la
T_{\a}^{\sg}\mid
\a\in\rd,\; \sg\in\cup_{J\in\supp(S)}(\tau_{_J} +2\Lam)\ra\\
&=&\la T_{\a}^{\sg}\mid \a\in\rd,\; \sg\in\tau_{_J}
+2\Lambda,\;J\in\supp(S)\ra\\ &=&\la
T_{\a}^{\tau_{_J}+\sg}\mid\a\in\rd,\;
\sg\in2\Lambda,\;J\in\supp(S)\ra\\ &=&\la
T_{\a}^{\tau_{_J}}T_{\a}^{\sg}T^{\tau_{_J}}_{\sg}\mid\a\in\rd,\;
\sg\in2\Lambda,\;J\in\supp(S)\ra.
\end{eqnarray*}
We get from  Lemma \ref{reform2-sim} and the facts that
$2\Lam\subseteq S$ and $\tau_{_J}\in\s$ for $J\in\supp(S)$, that
$$\begin{tabular}{c}
  ${\mathcal H}=\la T_{\a}^{\sg},~~
T_{\a}^{\tau_{_J}},\;T^{\tau_{_J}}_{\sg}\mid\a\in\rd,
\;\sg\in2\Lambda,\;J\in\supp(S)\ra$.
\end{tabular}$$
Now we show that each generator of ${\mathcal H}$ of the form
$T_{\a}^{\sg}$, $T_{\a}^{\tau_{_J}}$, $T^{\tau_{_J}}_{\sg}$
belongs to $T$. Let $\a=\sum_{i=1}^{\ell}m_{i}\a_{i}\in\rd$,
$\sg=\Sigma_{s=1}^{\nu}2n_s\sg_{s}\in 2\Lambda$, $n_s\in\bbbz$ and
$J\in\supp(S)$.
 Then it follows from Lemma \ref{formul-T},
 Lemma \ref{nice-lemma}, the definition of   $z_{_J}$ and the fact that
 $c_{r,s}^2\in\la z_{_{\{r,s\}}}\ra\subseteq T$,
 for all $1 \leq r,s\leq\nu$ that
 \begin{eqnarray*}
  T_{\tau_{_J}}^{\sg}&=&T_{\sum_{r\in
J}\sg_r}^{\Sigma_{s=1}^{\nu}2n_s\sg_{s}}=\prod_{r\in
J}\prod_{s=1}^{\nu}(c_{r,s}^{2})^{n_s}\in T
\end{eqnarray*}
and
\begin{eqnarray*}
T_{\a}^{\sg}&=&\prod_{r=1}^{\nu}\prod_{i=1}^{\ell}(t_{i,r})^
 {2m_{i}n_{r}}\prod_{1\leq r<s\leq\nu}
(c_{s,r}^{2})^{2n_{r}n_{s}}\in T.
 \end{eqnarray*}
Finally,
\begin{eqnarray*}
T_{\a}^{\tau_{_J}}&=&T_{\a}^{\sum_{r\in J}\sg_r}=\prod_{r\in
J}\prod_{i=1}^{\ell}t_{i,r}^
 {m_{i}}\prod_{\{r,s\in J\mid
r<s\}}c_{s,r}= \prod_{r\in J}\prod_{i=1}^{\ell}t_{i,r}^
 {m_{i}}z^{-1}_{_{J}}\in T.
 \end{eqnarray*}
From steps (1)-(2), the result follows.\qed

\begin{cor}\label{latice}
 If  $\ell\geq2$ or $S$ is a
lattice, then $$\mathcal H=\la t_{i,r},\;c_{r,s}\mid 1\leq
i\leq\ell,\;1\leq r\leq s\leq\nu\ra.$$
 \end{cor}
\proof This an immediate consequence of Proposition \ref{gen-H}
and the fact $F(S)=C$.\qed

The remaining results of this section are new only for type $A_1$.
In fact for types different from $A_1$, one can find essentially
the same results in \cite{MS} and \cite{A4}. However, for
completeness we provide a short proof of them, where the proofs
now are easy consequences of our results in Section
\ref{hyperbolic}.

\begin{pro}\label{propo-impor}

(i) $\w=\wdot\ltimes {\mathcal H}$.

(ii) If $\ell=1$, then $\w=\la w_{\a_1}, t_{1,r},z_{_J}\mid r\in
  J_\nu,\;J\sub J_\nu\ra.$

(iii) If $\ell>1$, then $\w=\la w_{\a_i}, t_{i,r},c_{r,s}\mid i\in
J_\ell,\;r\in
  J_\nu,\;1\leq r<s\leq\nu\ra.$

(iv) ${\mathcal H}$ is a torsion free group.

(v) ${\mathcal H}$ is a two-step nilpotent group.

 (vi) $Z(\w)=Z({\mathcal H})=F(S).$

(vii) $F(S)$  is a free abelian group of rank ${\nu(\nu-1)}/{2}$.
 \end{pro}

\proof (i) is an immediate consequence of Proposition \ref{xi-10}
and Lemma \ref{reform2-sim}. From (i), Corollary \ref{latice},
Proposition \ref{gen-H} and the fact that $\dot\w$ is generated by
$w_{\a_i},\ldots,w_{\a_\ell}$, it follows that (ii) and (iii)
hold. (iv) and (v) follow from  Lemma \ref{xi-4} and Corollary
\ref{xi-9}. From (\ref{def-F(S)}) and Proposition \ref{gen-H} we
have $F(S) \subseteq Z(\w)$. So to prove (vi) it is enough to show
that $Z(\w)\sub Z({\mathcal H})\sub F(S)$. Let $w\in Z(\w)$. By
(i), $w=\dot{w}h$ for some $\dot{w}\in \dot{\w}$ and $h\in
{\mathcal H}$. Since for all $(i,r)\in J_\ell\times J_\nu$,
$t_{i,r}\in {\mathcal H}$ we have from Lemma \ref{new-4} that
$$1=wt_{i,r}^{-1}w^{-1}t_{i,r}=wT_{-\a_{i}}^{\sg_{r}}w^{-1}T_{\a_{i}}^{\sg_{r}}=
T_{-w(\a_{i})}^{\sg_{r}}T_{\a_{i}}^{\sg_{r}}=T_{\a_{i}-w(\a_{i})}^{\sg_{r}}.$$
This gives $\a_i-w(\a_i)\in\v^0$. Since $w=\dot{w} h$ and
$h(\a_i)=\a_i$, $\mod\;\v^0$, it follows that $\dot{w}(\a_i)=\a_i$
for all $i\in J_\ell$. Thus $\dot{w}=1$ and so $w=h\in {\mathcal
H}\cap Z(\w).$ This gives $Z(\w)\subseteq Z({\mathcal H})$. Next
let $h\in Z({\mathcal H})$. From Proposition \ref{gen-H} and
(\ref{commutator-sim}), it follows that
$$h=z\prod_{r=1}^{\nu}\prod_{i=1}^{\ell}t_{i,r}^{m_{i,r}},\qquad(z\in
F(S),\;m_{i,r}\in\bbbz).$$ Then by (v), (\ref{bicomm}) and
(\ref{comu-sim}) we have that
\begin{eqnarray*}
1=[h,t_{j,s}]=\prod_{r=1}^{\nu}\prod_{i=1}^{\ell}[t_{i,r}^{m_{i,r}},t_{j,s}]
&=&\prod_{r=1}^{\nu}\prod_{i=1}^{\ell}
 c_{r,s}^{m_{i,r}(\a_{i},\a_{j})}=\prod_{r=1}^{\nu}
 c_{r,s}^{(\sum_{i=1}^{\ell}m_{i,r}\a_{i},\a_{j})}.
\end{eqnarray*}
Therefore from (\ref{free-abelian}), it follows that
$(\sum_{i=1}^{\ell}m_{i,r}\a_{i},\;\a_{j})=0$ for all $j\in
J_\ell$. Since the form on $\v$ restricted to
$\vd=\sum_{i=1}^{\ell}\bbbr\a_{i}$ is positive definite we get
$m_{i,r}=0$ for all $i\in J_\ell$ and $r\in J_\nu$. Then $h=z\in
F(S)$  and so (vi) holds. By (\ref{commutator-sim}) and the fact
that $c_{r,s}^{2}=[t_{1,r},t_{1,s}]$, we see that the group in the
statement is squeezed
 between two groups $\la c_{r,s}^2:1\leq
r<s\leq\nu\ra$ and $C$. Since $C$ is free abelian on generators
$c_{r,s}$, $1\leq r<s\leq\nu$, then (vii) follows. \qed

 The following
important type-dependent result gives explicitly the center
$Z(\w)=Z({\mathcal H})$ in terms of the generators $c_{r,s}$.
\begin{cor}\label{centers-1}
(i) If $X_\ell=A_1$, then $$Z({\mathcal H})=\la  c_{r,s}^2,\;
z_{_J}\mid 1\leq r< s\leq\nu,\;J\in\supp(S)\ra.$$ In particular if
$S$ is a lattice, then $Z({\mathcal H})=\la c_{r,s}\mid1\leq
r<s\leq\nu\ra.$

(ii) If $X_\ell=A_\ell(\ell\geq2), \;D_\ell$ or $E_\ell$, then
$Z({\mathcal H})=\la c_{r,s}\mid1\leq r<s\leq\nu\ra.$
\end{cor}

\proof Both (i) and (ii) follow immediately from (\ref{def-zJ})
and Proposition \ref{propo-impor}.\qed

\begin{pro}\label{w-form}

(i) If $X_\ell=A_1$, then each element $w$ in $\w$ has a unique
expression in the form
\begin{equation}\label{form=1}
  w=w(n,m_{r},z):=w_{\a_1}^n\prod_{r=1}^{\nu}t_{1,r}^{m_{r}}z\quad(n\in\{0,1\},\;m_{r}\in\bbbz,\;z\in F(S)).
  \end{equation}

(ii) If $X_\ell=A_\ell(\ell\geq2), \;D_\ell$ or $E_\ell$, then
each element $w$ in $\w$ has a unique expression in the form
\begin{equation}\label{form>1}
  w=w(\dot w,m_{i,r},m_{r,s} ):=\dot w\prod_{r=1}^{\nu}\prod_{i=1}^{\ell}t_{i,r}^{m_{i,r}}
  \prod_{1\leq r<s\leq\nu}c_{r,s}^{m_{r,s}},
  \end{equation}
  where $\dot w\in\dot\w,$ and $m_{i,r},\;m_{r,s}\in\bbbz$.
\end{pro}

 \proof (i) First we can express each element $w\in\w$ in terms
of generators given in Proposition \ref{propo-impor}(ii). Next we
can reorder the appearance of generators in any such expression
using (\ref{commutator-sim}), Proposition \ref{propo-impor}(vi),
Corollary \ref{centers-1}(i), Lemmas \ref{wtw-sim}  and the fact
that $w_{\a_1}^2=1$. Now to complete the proof it is enough to
show that the expression of $w$ in the form (\ref{form=1}) is
unique. Let $w(n',m_{r}',z')$
 be another expression of
$w$ in the form (\ref{form=1}). Then from  Proposition
\ref{propo-impor}(i) and Lemma \ref{xi-6}, it follows that $n=n'$
 $n_{r}=n_{r}'$ for all $r\in J_\nu$ and $z=z'$.

 (ii) Let $w\in\w$.  By parts (iii) and (vi) of  Proposition
\ref{propo-impor}, Corollary \ref{centers-1}(ii), Lemma
\ref{wtw-sim} and
 (\ref{commutator-sim}), $w$ can be written in the form
(\ref{form>1}).  Let $w(\dot{w}',n_{i,r}',m_{r,s}')$
 be another expression of
$w$ in the form (\ref{form>1}). Then  from  Proposition
\ref{propo-impor}(i) and Lemma \ref{xi-6}, it follows that $\dot
w=\dot w'$ and $n_{i,r}=n_{i,r}'$, $(i,r)\in J_\ell\times J_\nu$
and $m_{r,s}=m_{r,s}'$ for all $r,s\in J_\nu$. \qed

\vspace{5mm}
\section{\bf A  PRESENTATION FOR HEISENBERG-LIKE GROUP }
\label{presentation-sim}
 \setcounter{equation}{0}
We keep all the notations as in Section \ref{heisenberg-sim}. In
particular, $R$ is a simply laced extended affine root system and
$\mathcal H$ is its Heisenberg-like group.

We recall from Proposition \ref{propo-impor} that $F(S)=Z(\mathcal
H).$ For $1\leq r<s\leq\nu$, we define
\begin{eqnarray}\label{min-2-sim}
n(r,s)=\mbox{min}\{n\in\mathbb N: c_{r,s}^n\in F(S)\}.
\end{eqnarray}

\begin{rem}\label{rem1}
It is important to notice that we may consider $C=\la
c_{r,s}:1\leq r<s\leq\nu\ra$ as an abstract free abelian group
(see Lemma \ref{xi-7}(ii)) and $F(S)$ as a subgroup whose
definition depends only on the semilattice $S$. It follows that
the integers $n(r,s)$ are uniquely determined by $S$ (and so by
$R$). \end{rem}

We note from (\ref{def-zJ}) that if $1\leq r<s\leq\nu$, then
depending on either $\{r,s\}$ is in $\supp(S)$ or not we have
$z_{_{\{r,s\}}}=c_{r,s}$ or $z_{_{\{r,s\}}}=c_{r,s}^2$. Thus we
have the following lemma.

\begin{lem}\label{nrs-1-sim}
(i) $n(r,s)\in\{1,2\}$, $1\leq r<s\leq \nu.$

(ii) If $\{r,s\}\in\supp(S)$ for all $1\leq r<s\leq\nu$, then
$n(r,s)=1$ for all such $r,s$ and  $$F(S)=\la c_{r,s}~|~1\leq
r<s\leq\nu\ra.$$ In particular, this holds if $S$ is a lattice.

(iii) If $\supp(S)
 \subseteq\big\{\emp,\;\{r,s\}~~\mid~1<r\leq s\leq\nu\}$, then
 $$F(S)=\la c_{r,s}^{n(r,s)}~|~1\leq r<s\leq\nu\ra$$
 where
 \begin{eqnarray}\label{formula-1-sim}
n(r,s)=\left\{\begin{array}{ll}
 1, &\hbox{if}\; \{r,s\}\in\supp(S),\vspace{2mm} \\
2, &\hbox{if}\; \{r,s\}\not\in\supp(S).\\
 \end{array}\right.
 \end{eqnarray}
\end{lem}

Let $A=(a_{i,j})_{1\leq i,j\leq\ell}$ be the Cartan matrix of type
$X_\ell$, that is $$a_{i,j}=(\a_{i},\a_{j}),\quad i,j\in J_\ell.$$
\begin{lem}\label{div-1-sim}
$n(r,s)\mid a_{i,j},\quad i,j\in J_\ell$,\;\; $1\leq r<s\leq\nu$
\end{lem}
\proof By (\ref{commutator-sim}),
$c_{r,s}^{a_{i,j}}=c_{r,s}^{(\a_i,\a_j)}\in F(S)$ for all $r,s\in
J_\nu$ and $i,j\in J_\ell$. Now the result follows by the way
$n(r,s)$ is defined. \qed

By Lemma \ref{div-1-sim}, we have $a_{i,j}n(r,s)^{-1}\in\bbbz$,
$i,j\in J_\ell$, $1\leq r<s\leq\nu$.  So from (\ref{comu-sim}) we
have
\begin{equation}\label{comm-3-sim}
\begin{array}{c}
[t_{i,r},t_{j,s}]=(c_{r,s}^{n(r,s)})^{a_{i,j}n(r,s)^{-1}}.
\end{array}
\end{equation}
Note that the integer $n(r,s)$ appears only in type $A_1$ as in
other types $n(r,s)=1$.

\begin{thm}\label{presen-1-sim}
Let  $R=R(X_\ell, S)$ be a simply laced extended affine root
system  of type $X_\ell$ and nullity $\nu$  with Heisenberg-like
group $\mathcal H$. Let  $A=(a_{i,j})$ be the Cartan matrix of
type $X_\ell$, $n(r,s)$'s be the unique integers defined by
(\ref{min-2-sim}) and $m=\ind(S)$. If
\begin{equation}\label{000} F(S) \hbox{ is generated by elements
}c_{r,s}^{n(r,s)},\;\; 1\leq r<s\leq\nu,\end{equation} then
${\mathcal H}$ is isomorphic to the group $\widehat{\mathcal H}$
defined by generators
\begin{equation}\label{absh-sim}
\left\{
\begin{array}{ll}
y_{i,r}& 1\leq i\leq\ell,\;\;1\leq r\leq\nu,\\ z_{r,s}&1\leq
r<s\leq\nu,
\end{array}\right.
\end{equation} and
relations
\begin{equation}\label{relations-sim}
{\mathcal R}_{\widehat{\mathcal H}}:=\left\{\begin{array}{ll}
 [z_{r,s},z_{r',s'}],\vspace{2mm}\\

[y_{i,r},z_{r',s'}],\;[y_{i,r},y_{j,r}],\vspace{2mm}\\

[y_{i,r},y_{j,s}]=z_{r,s}^{a_{i,j}n(r,s)^{-1}},\;\;r<s,
\end{array}\right.
\end{equation}
where if $\ell>1$ or $\ell=1$ and $\nu\leq 3$, the condition
(\ref{000}) is automatically satisfied. Moreover, if $\ell
>1$, $n(r,s)=1$ for all $r,s$ and if $\ell=1$ and $\nu\leq 3$, $n(r,s)$'s are
given by the following table:

\begin{equation}\label{tab-100}
\begin{tabular}{c|c|c|c|c}
$\nu$&$m$ & $n(1,2)$ & $n(1,3)$ & $n(2,3)$ \\
  \hline
 0  & 0  & - & - & - \\
   \hline
 1  & 1  & - & - & - \\
     \hline
 2  & 2  & 2 & - & - \\
    & 3  & 1 & - & -\\
   \hline
 3  &  3 & 2 & 2 & 2 \\
    &  4 & 2 & 2 &1 \\
    &  5 & 2 & 1 & 1 \\
    & 6  & 1 & 1 & 1\\
    & 7  & 1 & 1 & 1\\
    \hline
  \end{tabular}
  \end{equation}
\end{thm}
\vspace{2mm}

 \proof By  Propositions  \ref{gen-H}, \ref{propo-impor}(vi),
(\ref{def-F(S)}) and assumption (\ref{000}), we have
  $${\mathcal H}=\la t_{i,r}\mid (i,r)\in J_\ell\times J_\nu\ra F(S)=\la
t_{i,r},c_{r,s}^{n(r,s)}\mid 1\leq i\leq\ell,\;1\leq
r<s\leq\nu\ra.$$ From (\ref{comm-3-sim}) and the fact that
$c_{r,s}^{n(r,s)}\in Z({\mathcal H})$, for $1\leq r<s\leq\nu$, it
is clear that there exists a unique epimorphism
$\varphi:\widehat{\mathcal H}\longrightarrow {\mathcal H}$ such
that $\varphi(y_{i,r})=t_{i,r}$ and $\varphi(z_{r,s})=
c_{r,s}^{n(r,s)}.$ We now prove that $\varphi$ is a monomorphism.
Let $\hat{h}\in\widehat{\mathcal H}$ and $\varphi(\hat{h})=1$. By
(\ref{relations-sim}), $\hat{h}$ can be written  as
\begin{equation*}
\hat{h}=\prod_{r=1}^{\nu}\prod_{i=1}^{\ell}y_{i,r}^{m_{i,r}}\prod_{1\leq
r<s\leq\nu} z_{r,s}^{n_{r,s}}\qquad(m_{i,r}, n_{r,s}\in\bbbz).
\end{equation*}
Then
$$1=\varphi(\hat{h})=\prod_{r=1}^{\nu}\prod_{i=1}^{\ell}t_{i,r}^{m_{i,r}}\prod_{1\leq
r<s\leq\nu} c_{r,s}^{n(r,s)n_{r,s}}.$$ Now it follows from
Proposition \ref{w-form} and  (\ref{free-abelian}) that
$m_{i,r}=n_{r,s}=0$, for all $i,r,s$ and so $\hat{h}=1$.

Next let $\ell>1$. By [AABGP, Proposition II.4.2] the involved
semilattice $S$ in the structure of $R$ is a lattice and so by
Lemma \ref{nrs-1-sim}(ii), $n(r,s)=1$ for all $r,s$. Finally, let
$\ell=1$. According to [AABGP, Proposition II.4.2], any extended
affine root system of type $A_1$ and nullity $\leq 3$ is
isomorphic to an extended affine root system of the form
$R=R(A_1,S)$  where $\supp(S)$ is given in Table \ref{tab-1}.
The result now follows immediately from this table and Lemma
\ref{nrs-1-sim}(ii)-(iii). \qed

\vspace{5mm}
\section{\bf A PRESENTATION FOR EXTENDED AFFINE WEYL GROUPS}\label{presentation-W-sim}
 \setcounter{equation}{0}
We keep the same notation as in the previous sections. As before
$R=R(X_\ell,S)$ is a simply laced extended affine root system of
nullity $\nu$,  $\w$ is its extended affine Weyl group and and
$\mathcal H$ is its Heisenberg-like group. Using the Coxeter
presentation for the finite Weyl group $\dot{\w}$, Theorem
\ref{presen-1-sim} and the semidirect product $\dot{\w}\ltimes
{\mathcal H}$, we obtain a finite presentation for $\w$. Let
$A=(a_{i,j})_{1\leq i,j\leq\ell}$ be the Cartan matrix of type
$X_\ell$.

We recall from \cite[ Proposition 3.13]{Ka} that $\dot\w$ is a
Coxeter group with generators $w_{\a_1},\ldots,w_{\a_\ell}$ and
relations
\begin{equation}\label{Coxe-rela}
   w_{\a_i}^2\andd (w_{\a_i}w_{\a_j})^{a_{i,j}^2+2}\quad (i\neq
  j).
 \end{equation}

\begin{thm}\label{presen-4-sim}
Let  $R=R(X_\ell, S)$ be a simply laced extended affine root
system  of type $X_\ell$ and nullity $\nu$  with extended affine
Weyl group $\w$. Let  $A=(a_{i,j})$ be the Cartan matrix of type
$X_\ell$ and $n(r,s)$'s be the unique integers defined by
(\ref{min-2-sim}).  If
\begin{equation}\label{000-w}
F(S) \hbox{ is generated by elements }c_{r,s}^{n(r,s)},\;\; 1\leq
r<s\leq\nu,\end{equation} then $\w$ is isomorphic to the group
$\widehat{\w}$ defined by generators
\begin{equation*}
 \left\{\begin{array}{ll}
x_i,&1\leq i\leq\ell, \\
 y_{i,r},& 1\leq i\leq\ell,\; 1\leq r\leq\nu,\\
 z_{r,s},&1\leq r<s\leq\nu\\
\end{array}\right.
\end{equation*} and
relations
\begin{equation*}
{\mathcal R}_{\widehat{\w}}:=\left\{\begin{array}{l}

x_i^2,\;(x_ix_j)^{a_{i,j}^2+2},\quad (i\neq
  j),\vspace{2mm}\\

x_{i}y_{j,r}x_i=y_{j,r}y_{i,r}^{-a_{i,j}},\vspace{2mm}\\

[z_{r,s},z_{r',s'}],\;[y_{i,r},z_{r',s'}],\;[y_{i,r},y_{j,r}],\vspace{2mm}\\

[y_{i,r},y_{j,s}]=z_{r,s}^{a_{i,j}n(r,s)^{-1}},\;1\leq r<s\leq\nu.
\end{array}\right.
\end{equation*}
Moreover if $\ell>1$ then $n(r,s)=1$ for all $r,s$, (in
particular, the assumption (\ref{000-w}) holds). Furthermore, if
$\ell=1$ and $\nu\leq 3$ then the assumption (\ref{000-w}) is
automatically satisfied and the relations $\rr_{\widehat{\w}}$
reduces to the relations
\begin{equation*}
{\mathcal R}_{\widehat{\w}}:=\left\{\begin{array}{l}

x^2,\;\;xy_{r}x=y_{r}^{-1},\vspace{2mm}\\

[z_{r,s},z_{r',s'}],\;[y_{r},z_{r',s'}],\vspace{2mm}\\

[y_{r},y_{s}]=z_{r,s}^{2n(r,s)^{-1}},\;1\leq r<s\leq\nu,
\end{array}\right.
\end{equation*}  where $n(r,s)$'s are given explicitly
by (\ref{tab-100}) (depending on $m=\ind(S)$).
 \end{thm}

\proof From parts (ii), (iii) and (vi) of  Proposition
\ref{propo-impor},  (\ref{def-F(S)}) and assumption (\ref{000-w}),
it follows that
\begin{equation}
\begin{array}{ll}
 \w=\la w_{\a_i},\; t_{i,r}\mid (i,r)\in J_\ell\times J_\nu\ra
F(S)\vspace{2mm}\\
 \hspace{4mm}=\la w_{\a_i},\; t_{i,r},\;c_{r,s}^{n(r,s)}\mid 1\leq
i\leq\ell,\;1\leq r<s\leq\nu\ra. \\
\end{array}
 \end{equation}
 By (\ref{comm-3-sim}), Lemma \ref{wtw-sim}, (\ref{Coxe-rela}) and
the fact that $c_{r,s}^{n(r,s)}\in Z({\mathcal H})$, for $1\leq
r<s\leq\nu$, the assignment $x_i\longmapsto w_{\a_i}$,
$y_{i,r}\longmapsto t_{i,r}$ and $z_{r,s}\longmapsto
c^{n(r,s)}_{r,s}$  induces a unique epimorphism $\psi$ from
$\widehat{\w}$ onto $\w$. Also by Lemma (\ref{Coxe-rela}), the
restriction of $\psi$ to $\widehat{\dot\w}:=\la x_i\mid 1\leq
i\leq\ell\ra$ induces the isomorphism
 \begin{equation}\label{finite-case}
 \widehat{\dot\w}\cong^{^{\hspace{-2mm}\psi}}\dot{\w}.
 \end{equation}
 We now show that $\psi$ is injective.
Let $\psi(\hat{w})=1$, for some $\hat{w}\in\widehat{\w}$. From the
defining relations for $\widehat{\w}$, it is easy to see that
$\hat{w}$ can be written in the form
\begin{equation*}
  \hat{w}=\hat{\dot{w}}\prod_{r=1}^{\nu}\prod_{i=1}^{\ell}y_{i,r}^{m_{i,r}}
  \prod_{1\leq r<s\leq\nu}z_{r,s}^{n_{r,s}}
  \quad(\hat{\dot{w}}\in\widehat{\dot{\w}},\hspace{1mm}  n_{r,s},
  m_{i,r}\in\bbbz).
  \end{equation*}
 Then $$ 1=\psi(\hat{w})=\psi(\hat{\dot{w}})\prod_{r=1}^{\nu}
  \prod_{i=1}^{\ell}t_{i,r}^{m_{i,r}}
  \prod_{1\leq r<s\leq\nu}c_{r,s}^{n(r,s)n_{r,s}}.$$
  Therefore from (\ref{finite-case}) and Propositions \ref{w-form}, it
follows that  $m_{i,r}=0$, $n_{r,s}=0$ for all $i,r,s$  and
$\hat{\dot w}=1$. Thus $\hat{w}=1$ and so $\w\cong\widehat{\w}$.
Now an argument similar to the last paragraph of the proof of
Theorem \ref{presen-1-sim} complets the proof. \qed

 We close this section with the
following remark.

\begin{rem}\label{rem2}
In Section \ref{heisenberg-sim}, we fixed a hyperbolic extension
$(\vt,I)$ of $(\v,I)$, determined by extended affine root system
$R$ and then we defined the extended affine Weyl group $\w$ as a
subgroup of $\O(\vt,I)$. However, Remark \ref{rem1} and Theorem
\ref{presen-4-sim} show that the definition of $\w$ is independent
of the choice of this particular hyperbolic extension.
\end{rem}


\begin{thebibliography}{AABGP}
\bibitem[AABGP]{AABGP} {\it B. Allison, S. Azam, S. Berman,
Y. Gao, and A. Pianzola, Extended
 affine Lie algebras and their root systems}, Mem. Amer. Math. Soc. {\bf 603}
(1997), 1-122.

\bibitem[A1]{A1} S. Azam, {\it Nonreduced extended affine Weyl groups}, J.
Alg. {\bf 269} (2003), 508-527.

\bibitem[A2]{A2} S. Azam, {\it Extended affine root systems},  J. Lie
Theory 12 (2) (2002), 515-527.

\bibitem[A3]{A3} S. Azam, {\it A presentation for reduced extended affine Weyl groups},
Comm. Alg. {\bf 28} (2000), 465-488.

\bibitem[A4]{A4} S. Azam, {\it Extended affine Weyl groups}, J. Algebra
{\bf 214} (1999) 571-624.

\bibitem[A5]{A5} S. Azam, {\it Nonreduced extended affine root system of nullity
3}, Comm. Alg. {\bf 25} (1997), 3617-3654.


\bibitem[H] {H} J. E. Humphreys, {\it Reflection groups and Coxeter groups}
, Cambridge University Press, 1992.


\bibitem[Ka]{Ka} V.G. Kac,{\it Infinite dimensional Lie algebras},
Third edition, Combridge University Press, Cambridge, 1990.

\bibitem[Kr]{Kr} Y. Krylyuk, {\it On the structure of quasi-simple Lie
algebras and their automorphism groups}, Ph.D Theses, University
of Saskatchewan, 1995.


\bibitem[MS]{MS} R. V. Moody and Z. Shi, {\it Toroidal Weyl groups}, Nova~J.~Algebra
~Geom. {\bf 11} (1992), 317-337.



\bibitem[S]{S} K. Saito, {\it Extended affine root systems 1(Coxeter
transformations)}, RIMS., Kyoto Univ. {\bf 21} (1985), 75-179.

\bibitem[ST]{ST} K. Saito and T. Takebayashi, {\it Extended affine
root systems III (elliptic Weyl groups)},
Publ.~Res.~Inst.~Math.~Sci. {\bf 33} (1997), 301-329.

\bibitem[T]{T} T. Takebayashi,
{\it Defining relations of Weyl groups for extended affine root
systems $A^{(1,1)}_\ell$, $B^{(1,1)}_\ell$, $C^{(1,1)}_\ell$,
$D^{(1,1)}_\ell$}, J. Alg. {\bf 168} (1994), 810--827.

\end{thebibliography}
\end{document}